\documentclass[12pt]{article}
\usepackage{amssymb}
\newtheorem{theorem}{Theorem}
\newtheorem{lemma}{Lemma}

\begin{document}
\title{On reduction integer programs to knapsack problem.}\author{S.I. Veselov}
\date{N.I. Lobachevsky State University of Nizhni Novgorod}
\maketitle

{\bf Abstract.}
{\small Let $A$ be an integral nonnegative $m\times n$ matrix, $b$ be an integral nonnegative vector. It is suggested
new method for reduction of integer program $\max \{cx|\  Ax=b,\ x\ge 0,\ x\in\mathbf{Z}^n\}$ to knapsack problem
$\max \{c'x|\  fAx=fb,\ x\ge 0,\ x\in\mathbf{Z}^n\}$.  
}
\\\bigskip

{\bf Introduction.} Let $A$ be an integral nonnegative $m\times n$ matrix, $b=(b_1,..,b_m)^T$ be an integral nonnegative vector. 
By definition, put $M(A,b)=\{x| Ax= b,\ x - {\rm integral\  and\ nonnegative}\}$. Denote by $V(A,b)$ the vertex set of 
convex hull of $M(A,b)$.
Integral vector $f$ is called {\it aggregating } for $M(A,b)$ if $M(A,b)=M(f^TA,f^Tb).$

In [1], [2] several aggregating vectors have been described. The bibliography see in [1]. In [2],[3] the lower estimates for coordinates of the aggregating vector have been given. 

It is clear that if $f$ is aggregating vector for $M(A,b)$ then $\max \{c^Tx|\  Ax=b,\ x\ge 0,\ x\in\mathbf{Z}^n\}=\max \{cx|\  f^TAx=f^Tb,\ x\ge 0,\ x\in\mathbf{Z}^n\}$.  

In this paper it is suggested
new approach for reduction of integer program $\max \{c^Tx|\  Ax=b,\ x\ge 0,\ x\in\mathbf{Z}^n\}$ to knapsack problem
like $\max \{c'x|\  f^TAx=f^Tb,\ x\ge 0,\ x\in\mathbf{Z}^n\}$ .
\\\medskip

{\bf Main result.}
Integral vector $f$ is called {\it V-aggregating} for $M(A,b)$ if $V(A,b)\subseteq V(f^TA,f^Tb).$
Put by definition $f=(1,b_1+1,..,\prod\limits_{i=1}^{m-1}(b_i+1))^T, e=(1,1,..,1)^T$

\begin{lemma}\label{l1} $b\in V(f,f^Tb). $\end{lemma}
Proof.
It is sufficient to note that linear function $h(t)=e^Tt$ achieves its minimum over 
$M(f,f^Tb)$ at the unique vector $t=b$.

\begin{theorem} If $M(A,b)\ne \emptyset$ then $V(A,b)\subseteq V(f^TA,f^Tb).$
 \end{theorem}

Proof. Let us show that  $x^0\notin V(f^TA,f^Tb)$ implies $x^0\notin V(A,b)$
Let $x^1,..x^r\in M(f^TA,f^Tb)$ 

$$
x=\sum\limits_{i=1}^r\lambda_ix^i,\,\sum\limits_{i=1}^r\lambda_i=1,\,\lambda_i\ge
0\,(i=1,..,r)\eqno(1)$$ Then
$$
Ax^0=\sum\limits_{i=1}^r\lambda_iAx^i,\,\sum\limits_{i=1}^r\lambda_i=1,\,\lambda_i\ge
0\,(i=1,..,r)\eqno(2)$$ Since $Ax^i\in M(f,f^Tb) (i=1,..,r)$  it follows from Lemma {\ref{l1}}
that $Ax^i=b$, hence,  $x_i\in M(A,b)(i=1,..,r)$.
From (1) it follows that  $x^0\notin V(A,b)$.

\paragraph{Reduction to knapsack. }

Let $L\ge \min\limits_{x\in M(A,b)}cx $,  $k$ satisfies the inequality $c+ke^TA\ge
0$ and $H=L+ke^Tb+1$.

\begin{theorem}
If $M\ne \emptyset$ and linear function $c^Tx+He^TAx$ achieves its minimum over $M(f^TAx,f^Tb)$ at $x^0$
then $x^0\in M(A,b)$ and $c^Tx^0=\min\limits_{x\in M(A,b)}c^Tx.$
\end{theorem}
Proof. Suppose $Ax^0\ne b$. It follows from Lemma \ref{l1} that
$$c^Tx^0+He^TAx^0\ge -ke^TAx^0+He^TAx^0\ge (H-k)(e^Tb+1)=k(e^Tb)^2+(L+k+1)e^Tb+L+1.$$ On the other hand, there exists $x'\in M(A,b)$ such that $c^Tx'\le L$, hence $$c^Tx'+He^TAx'\le L+He^Tb=k(e^Tb)^2+(L+k+1)e^Tb+L.$$

\paragraph{Low bound}

\begin{theorem} If $a$ is nonnegative intergal vector, $V(a^T,a_0)\subseteq V(A,b),\ x^0=(x^0_1,..,x^0_n)\in V(A,b)$ then
$a_0\ge \prod\limits_{i=1}^n(x_i^0+1)-1$ \end{theorem}
Proof. Let us consider $t^1,t^2$ such that
$$t^1\ne t^2, 0\le t^1\le x^0, 0\le t^2\le x^0, a^Tt^1=a^Tt^2.
$$ Since $x^1=x^0-t^1+t^2, x^2=x^0+t^1-t^2\in M(a^T,a_0)$ and $x^0=\frac{1}{2}(x^1+x^2)$. But this contradicts the 
fact that $x^0\in V(A,b).$
 So the linear function $a^Tt$  over $0\le
t\le x^0$ takes $\prod\limits_{i=1}^n(x_i^0+1)$ distinct nonnegative integers not exceeding $a_0$.

\end{document}